\documentstyle[12pt]{article}

\begin{document}
\newcommand{\nt}{\noindent}

\begin{center}

{\bf Another way to enumerate rational curves with
torus actions}

\bigskip

Aaron Bertram  

\end{center}

\medskip

\nt {\bf 1. Introduction.} In 1991 Candelas, de la Ossa, 
Green and Parkes used the ``string-theoretic'' principle of mirror
symmetry to predict the numbers of rational curves of any degree on a
quintic three-fold  
\cite{COGP}. This prediction took mathematicians by surprise, as the
best results at the time only counted rational curves
of degree three or less. Since then, exciting new developments 
have led to mathematical proofs of these predictions and many other
related conjectures coming from string theory. While this does not address
the deeper problem of constructing a mathematical foundation for string
theory, it does represent a major advance in field of 
enumerative algebraic geometry.   

\medskip  

An enumerative question is usually
interpreted in terms of intersection theory on a moduli
space. Moduli spaces of stable maps were introduced by Kontsevich and
Manin \cite{KM} on which to calculate Gromov-Witten invariants, including
the expected numbers of rational curves on a quintic three-fold. Like the
spaces of stable pointed curves (which are stable maps to a point), the 
boundary of a stable map space has a self-similar property that has been
exploited in many of the important recent results in enumerative 
geometry, including the associativity of the quantum cohomology rings
\cite{KM,RT} and the reconstruction
theorems \cite{Du,KM} for genus zero invariants.   
 
\medskip

Any torus action on the target space is inherited by all spaces
of stable maps. 
Kontsevich applied the Bott residue theorem for such torus
actions in his {\it Enumeration of rational curves with torus actions}
\cite{Ko} to the problem of computing Gromov-Witten invariants of
rational curves, and although his approach would in principle compute the 
numbers of rational curves of all degrees on the quintic, the
combinatorics  becomes unmanageable after degree four. Givental seems
to have been the first to fully appreciate the significance of
the more subtle circle action induced from the {\bf domain} of a 
``parametrized'' stable
map. An application of techniques from equivariant cohomology to both
actions led   to Givental's proofs
\cite{BDPP,Gi,Gi2,Gi3,Pa} of the mirror conjecture for Fano and
Calabi-Yau complete intersections in toric varieties, as well as the
subsequent proofs and generalizations of Kim \cite{Ki} and Lian-Liu-Yau
\cite{LLY,LLY2}.   In particular, these verify the predicted numbers of 
rational curves on the quintic three-fold. 

\medskip

In this paper we will produce a new and direct computation of 
the one-point Gromov-Witten invariants for rational curves on Fano and
Calabi-Yau complete interesections in complex projective space. This proof
is self-contained and much
simpler than previous proofs. The idea is to use the self-similar
properties of the boundary of stable map space within the context of
equivariant cohomology to find  a self-similar decomposition of the
relevant  equivariant ``virtual'' classes. This new idea allows us to
dispense entirely with torus actions inherited from the target manifold
and instead focus on the more important circle action inherited from the
parametrization of the genus zero curve. This proof has the additional
advantage in that it gives a computation that does not rely on ``already
knowing the answer,'' as is the case in all previous proofs of the mirror
conjecture. Besides  being a psychological advantage, our proof thereby
generalizes without modification to give a relative version of the mirror
conjecture for projective bundles.

\medskip

The mirror conjecture is best expressed in terms of one-point invariants.
The space of stable genus-zero one-pointed maps of degree $d$ to ${\bf
P}^n$ is denoted by
$\overline M_{0,1}({\bf P}^n,d)$ and comes equipped with 
evaluation and projection maps:
$$e:\overline M_{0,1}({\bf P}^n,d) \rightarrow {\bf P}^n
\ \mbox{and}\ \pi:\overline M_{0,1}({\bf P}^n,d) \rightarrow 
\overline M_{0,0}({\bf P}^n,d)$$
where $\overline M_{0,0}({\bf P}^n,d)$ is the space of stable
zero-pointed maps. These spaces are always projective orbifolds and when
$d=1$ the maps are the two projections:
$$e:{\rm Fl}(1,2,n+1) \rightarrow {\bf P}^n \ \mbox{and}\ 
\pi:{\rm Fl}(1,2,n+1) \rightarrow {\rm G}(2,n+1)$$ 
from the partial flag manifold.

\medskip

If $S \subset {\bf P}^n$ is a transverse zero section of a vector
bundle $E$ which is generated by global sections (e.g. a complete
intersection), then
$S$ defines a ``virtual'' class:
$$[S]_d := \pi^*c_r(\pi_*e^*E)$$
on $\overline M_{0,1}({\bf P}^n,d)$, where $r$ is the rank of
$\pi_*e^*E$. This measures the expected constraint imposed by 
requiring the stable map to land in $S$. Even simpler, given a variety $V
\subset {\bf P}^n$ with associated cohomology class $[V]$, the pull-back:
$$e^*[V]$$  
is the constraint imposed by requiring the marked point to land
in $V$. As opposed to the previous two ``enumerative'' classes, Witten's
cotangent class:
$$\psi := c_1(\omega_\pi)$$
(the first chern class of the relative dualizing
sheaf for the projection $\pi$) is of less enumerative significance,
but very useful for computations.

\medskip

The mirror conjecture for ${\bf P}^n$ (as formulated by
Givental and generalized to homogeoneous spaces by Kim \cite{Ki}) says
that if $S$ is a Fano or Calabi-Yau complete intersection in
${\bf P}^n$ of type
$(l_1,...,l_m)$, then
the Laurent polynomials:
$$e_*\left(\frac{[S]_d}{t(t-\psi)}\right) := t^{-2}e_*([S]_d) +
t^{-3}e_*([S]_d \cup \psi) + t^{-4}e_*([S]_d \cup \psi^2) + ...$$ 
with coefficients in H$^*({\bf P}^n,{\bf Q})$ are closely related to the
rational functions (also expressible as Laurent polynomials):
$$\phi_d(t,h) := \frac{\prod_{i=1}^m\prod_{k=0}^{dl_i}(l_ih +
kt)}{\prod_{k=1}^d(h + kt)^{n+1}}$$
where $h\in \mbox{H}^*({\bf P}^n,{\bf Z})$ is the hyperplane class.

\medskip

More precisely, the relationship depends on
the positivity of $S$ as follows:

\medskip

\nt {\bf  Fano of Index Two or More:} If $l_1 + ... + l_m <
n$, then:
$$e_*\left(\frac{[S]_d}{t(t-\psi)}\right) = \phi_d$$

\nt {\bf Fano of Index One:} If $l_1 + ... + l_m = n$,
then:
$$e_*\left(\frac{[S]_d}{t(t-\psi)}\right) = 
\sum_{r=0}^d \frac{(-\prod l_i!)^{r}\phi_{d-r}}{r!t^r}$$

\nt {\bf Calabi-Yau:} If $l_1 + ... + l_m = n+1$, then there exist
power series:
$$f(q) = \sum_{d=1}^\infty a_dq^d \ \ \mbox{and}\ \ g(q) =
\sum_{d=1}^\infty b_dq^d$$ 
with constant coefficients such that the two power series:
$$\Phi(q) =
\sum_{d=0}^\infty
\phi_d q^d \ \
\mbox{and}\ \  \Sigma(q) = [S] + \sum_{d=1}^\infty e_*\left(
\frac{[S]_d}{t(t-\psi)}\right)q^d$$  
(with Laurent polynomial coefficients) satisfy:
$$\Sigma(q) = e^{\frac ht f(q) + g(q)}\Phi(qe^{f(q)})$$ 

If we expand the right side, we get:
$$e^{\frac ht f(q) + g(q)}\Phi(qe^{f(q)}) = 
\sum_{d=0}^\infty \phi_d(q^de^{(d + \frac ht)f(q) + g(q)})$$
$$= \sum_d\sum_{0 \le d_1 < d_2 < ... < d_{r+1} = d}
\frac{\phi_{d_1} \prod_{i=1}^r (a_{d_{i+1} - d_i}(d_1 + \frac ht)
+ b_{d_{i+1} - d_i})}{r!}q^d$$
so that the mirror conjecture is equivalent to the existence of 
constants $a_d$ and $b_d$ for $d = 1,...,\infty$ such that:
$$e_*\left(\frac{[S]_d}{t(t-\psi)}\right) = \sum_{0 \le d_1 < d_2
< ... < d_{r+1} = d}
\frac{\phi_{d_1} \prod_{i=1}^r (a_{d_{i+1} - d_i}(d_1 + \frac ht)
+ b_{d_{i+1} - d_i})}{r!}$$ 

All the previous proofs proceed by finding sufficient
conditions shared by the power series $\Sigma(q)$ and the change of
variables of the power series $\Phi(q)$ to uniquely characterize them,
hence to  conclude that they are the same. It is in this sense that one
needs to ``know the answer'' beforehand. Our proof is completely
different. We will explicitly produce the coefficients of the power
series one coefficient at a time.  Out of the proof it becomes clear how
the  self-similar properties
of the coefficients of the $\Sigma(q)$ in the Calabi-Yau 
and Fano of index one cases are a reflection
of the self-similar properties of the boundary of the
moduli spaces of stable 
maps.

\bigskip

{\bf Acknowledgements:} This paper grew out of a seminar in the summer of
1998 co-organized with Holger Kley during which
Chris Hacon and Herb Clemens made many useful suggestions. The use of
MacPherson's graph construction in Lemma 4.4 was suggested by Bill Fulton.
Alexander Givental was very helpful when I reported on the Fano version
of the result at MSRI in December 1998, and the key step in the Calabi-Yau
proof arose out of a conversation with Ravi Vakil.
Rahul Pandharipande, Pavel Etinghof and especially Bumsig Kim
helped  relate my original Calabi-Yau formula to the mirror conjecture and
Bumsig Kim  and Bong Lian made some useful remarks on an
earlier version of the paper. Holger Kley very helpful with a
critical revision of the paper and I'd like to thank  
Michael Thaddeus for explaining the basics of
equivariant cohomology one sunny afternoon in Park City.

\medskip

{\bf Note:} A very recent preprint of Gathmann \cite{Ga}
explains how one might compute
rational Gromov-Witten invariants of hypersurfaces without any torus
actions. This seems promising, particularly in the general-type
case.  
 
\newpage

\nt {\bf 2. Genus Zero Gromov-Witten Invariants on Projective Space.} 

\medskip

\nt {\bf Definition:} If $(C;p_1,...,p_k)$ is a pointed genus-zero curve,
then $f:C \rightarrow {\bf P}^n$
is stable if:

\medskip

(a) $C$ has only nodes for singularities and the $p_i$ are nonsingular
points.

\medskip

(b) Every component of $C$ which is collapsed by $f$ has at least 
three distinguished points (a point is distinguished
if it is either a node or marked).

\medskip

\nt {\bf Existence Theorem (see \cite{KM},\cite{Al}):} The moduli problem
for  isomorphism classes of flat families of genus-zero 
$k$-pointed stable maps of degree $d$ to ${\bf
P}^n$ is represented by a projective orbifold. This moduli space is
denoted by
$\overline M_{0,k}({\bf P}^n,d)$
with universal curve ${\cal C} \cong \overline M_{0,k+1}({\bf P}^n,d)$
and maps:
$$ \overline M_{0,k}({\bf P}^n,d)\stackrel \pi\leftarrow{\cal C} \stackrel
e\rightarrow {\bf P}^n$$ 

For each $1 \le i \le k$, if 
$\sigma_i$ is the section of $\pi$ determined $p_i$, then
$e_i := e \circ \sigma_i$ is the evaluation of a stable map at the marked
point $p_i$.

\medskip

The relevant Chern classes on $\overline M_{0,k}({\bf P}^n,d)$
come from three sources. 

\medskip

$\bullet$ The pull-back classes $e_i^*(h^b)$ 

\medskip

$\bullet$ Witten's cotangent classes $\psi_i := c_1(N^*_{\sigma_i})$ for
$i = 1,...,k$. When $k = 1$, then $\psi := \psi_1$ is the relative
canonical class $c_1(\omega_\pi)$. 

\medskip

$\bullet$ The top chern class $c_r(\pi_*e^*(E))$ whenever 
$E$ is a vector bundle on ${\bf P}^n$ which is generated by its global
sections. If $S \subset {\bf P}^n$ is the zero locus of a general section
of $E$, then this chern class will be denoted by $[S]_d$ as before.

\bigskip

The classes $[S]_d$ are all pulled back from $\overline M_{0,0}({\bf
P}^n,d)$, leading to relations: 

\medskip

\nt {\bf Example:} 
If $S \subset {\bf P}^4$ is a quintic hypersurface, then the 
``physical'' number of rational curves of degree $d$ on $S$ is given by:
$$\int_{\overline M_{0,0}({\bf P}^4,d)} [S]_d = 
\frac 1{d}\int_{\overline M_{0,1}({\bf P}^4,d)} e^*(h) \cup [S]_d
= \frac 1d \int_{{\bf P}^4}e_*[S]_d \cup h$$
so that the cohomology classes $e_*[S]_d$ determine all the 
``physical'' numbers of rational curves on the quintic (we will mostly 
ignore the subtleties of multiple covers in this paper).

\medskip

The general $k$-point genus zero Gromov-Witten invariant of a
complete intersection $S\subset {\bf P}^n$ of type $(l_1,...,l_m)$ is
an intersection number of the form:
$$\int_{\overline M_{0,k}({\bf P}^n,d)} p(e_i^*(h),\psi_i) \cup 
[S]_d$$  where $p(x_i,y_i)$ is a polynomial (or power series) in $2k$
variables.

\medskip

Notice that in contrast to zero-point invariants, there are interesting
one-point invariants on {\bf any} complete intersection (because of the
$\psi$'s), including projective space itself. 

\medskip

\nt {\bf Definition:} Let $(C \supset {\bf P}^1;p_1,...,p_k)$ be a
genus-zero curve with $k$ marked points and a distinguished parametrized
component ${\bf P}^1 \subset C$. Then a map
$f:C \rightarrow {\bf P}^n$
is stable if:

\medskip

(a) $C$ has only nodes for singularities and the $p_i$ are nonsingular
points.

\medskip

(b) Every component of $C$ which is collapsed by $f$ is either 
distinguished or else has at
least three distinguished points. 

\bigskip

The existence theorem applies in the context of stable maps with a
parametrized component. Following Givental, we'll call these moduli spaces
the ``graph spaces'' and denote them by
$\overline N_{0,k}({\bf P}^n,d)$.
In addition to the maps of the theorem, the universal curve over the graph
space admits an evaluation map $\epsilon:{\cal C} \rightarrow {\bf P}^1$
leading to corresponding evaluation maps at the points. 
The zero-pointed graph space
is rational via an explicit birational morphism:
$$u: \overline N_{0,0}({\bf P}^n,d) \rightarrow 
{\bf P}^n_d := {\bf P}^{(n+1)(d+1) - 1}$$
defined as follows. The general point of $\overline N_{0,0}({\bf P}^n,d)$
is represented by a map $f:{\bf P}^1 \rightarrow {\bf P}^n$ of degree $d$.
Such a map is given by $n+1$ polynomials
of degree $d$
and thus represents a point in the projective space 
${\bf P}^n_d$. Boundary points of the graph space 
are represented by maps $f:C \rightarrow {\bf P}^n$ where $C$
has several components, one of which is ${\bf P}^1$. But if $q_1,...,q_k
\in {\bf P}^1$ are the nodes of $C$ on ${\bf P}^1$, and if the
curve $C_i$ growing  out of
$q_i$ is mapped to ${\bf P}^n$ with degree $d_i$, then $f$ maps
${\bf P}^1$ to ${\bf P}^n$ with degree $d - \sum d_i$. Letting 
this map be given by polynomials $p_i$, we get a
well-defined point of ${\bf P}^n_d$ by choosing linear
forms $q_i(x,y)$ dual to the points $q_i$ and sending the stable map $f$
to:
$$(\prod_{i=1}^k q_i(x,y)^{d_i}p_0(x,y):....:\prod_{i=1}^k
q_i(x,y)^{d_i}p_n(x,y))$$
(It is easy to exhibit $u$ as a morphism. See Jun Li's argument in
\cite{LLY})   There is an inclusion:
$$i: \overline M_{0,1}({\bf P}^n,d) \hookrightarrow \overline
N_{0,0}({\bf P}^n,d)$$
defined as follows. Given a stable map $f:C \rightarrow {\bf P}^n$ with one
marked point $p\in C$, construct a new curve $C' = C \cup {\bf P}^1$ by 
joining ${\bf P}^1$ and $C$ at $p \in C$ and 
$0 \in {\bf P}^1$. Extend $f$ to a map  $f':C' \rightarrow {\bf P}^n$ by
collapsing the ${\bf P}^1$ component to the image point $f(p)$. Then there
is a:

\medskip

\nt {\bf Basic Diagram:} 
$$\begin{array}{ccc} \overline N_{0,0}({\bf P}^n,d)  & \stackrel u
\rightarrow & {\bf P}^n_d\\ \\ 
i \uparrow \ \ & & j\uparrow \ \ \\ \\   
 \overline M_{0,1}({\bf P}^n,d) &  \stackrel e \rightarrow &
{\bf P}^n
\end{array}$$
where $j(a_0:...:a_n) = (a_0x^d:...:a_nx^d)$.

\bigskip

\nt {\bf 3. Equivariant Basics.} Let $X$ be a compact complex 
manifold (or orbifold) equipped with a ${\bf C}^*$ action. (We will
let $T = {\bf C}^*$.) Then:
$$ET := {\bf C}^{\infty + 1} - \{0\} \rightarrow {\bf CP}^\infty =: BT$$
is the ``universal'' principal ${\bf C}^*$ bundle, which we use to
construct:
$$\pi_X:X_T := X\times_{{\bf C}^*} ET \rightarrow  ET/{\bf C}^* = BT$$

Following Borel, we define the equivariant cohomology by setting:
$$H^*_T(X,{\bf Q}) := H^*(X_T,{\bf Q})$$
which is a module via $\pi_X^*$ over
$H^*(BT,{\bf Q}) \cong {\bf Q}[t]$. 
At one extreme, if the $T$ action were trivial, this would be
the tensor product $H^*(X,{\bf Q})\otimes_{\bf Q} {\bf Q}[t]$ while
at the other extreme if the  action were free, it would be the 
cohomology ring of the quotient $H^*(X/T,{\bf Q})$, which is torsion as a
module over ${\bf Q}[t]$. 

\medskip

A vector bundle $E$ over $X$ is linearized if it is
equipped with an action of $T$ which is lifted from the action on $X$ so
that:
$$E_T \rightarrow X_T$$
is a vector bundle which pulls back to $E \rightarrow X$ on each of the
fibers of $X_T$ over $BT$. As such, its chern classes represent elements of
the equivariant cohomology ring, and one defines:
$$c_i^T(E) := c_i(E_T) \in H^*_T(X,{\bf Q})$$

If $f:X \rightarrow X'$ is a $T$-equivariant morphism of compact complex
manifolds with $T$ actions, let $f:X_T \rightarrow X'_T$ also denote the 
induced map on these spaces (which commutes with the projections to $BT$).
The pull-back on equivariant cohomology is the ordinary
pull-back with respect to the induced map $f$:
$$f^*:H^*_T(X',{\bf Q}) \rightarrow H^*_T(X,{\bf Q})$$
and if $f$ is proper, then the equivariant
proper push-forward is defined in the same way. We may now state: 

\medskip

\nt {\bf The Atiyah-Bott Localization Theorem:} Let $X$ be a compact
complex manifold with a $T$-action, let
$F_1,...,F_n \subset X$ be the (necessarily smooth) connected components
of the fixed-point locus
and let $i_k:F_k \hookrightarrow X$ denote their embeddings. Each
normal bundle to $F_k$ is canonically linearized, its top equivariant
chern class
$\epsilon^T(N_{F_k/X})$ is invertible in $H^*(F_k,{\bf Q})[t,t^{-1}]$,
and every torsion-free element $c_T \in H^*_T(X,{\bf Q})$
uniquely  decomposes in $H^*_T(X,{\bf Q})\otimes {\bf Q}(t)$ as a sum of
contributions from fixed-point loci:
$$c_T = \sum_{k=1}^n (i_k)_* \frac{i_k^*(c_T)}{\epsilon^T(N_{F_k/X})}$$

\nt {\bf Remark:} The localization theorem as stated here and proved in
\cite{AB} only  applies to compact complex manifolds $X$. Generalizations
to the case where $X$ is an orbifold may be found in the
papers of Graber-Pandharipande \cite{GP} and Kresch \cite{Kr}. See
\cite{AP} for a general reference on equivariant cohomology. 

\medskip
 
\nt {\bf Corollary (Correspondence of Residues):} Suppose $f:X
\rightarrow Y$ is an equivariant map of compact complex manifolds 
(or orbifolds) with $T$ actions and suppose $i:F \hookrightarrow X$ and
$j:G \hookrightarrow Y$ are components of the fixed-point loci with the
property that $F$ is the unique fixed component to map to $G$. Then
equivariant cohomology  classes $c_T$ on $X$ satisfy:
$$(f|_F)_*\left(\frac{i^*c_T}{\epsilon^T(N_{F/X})}\right) = 
\frac{j^*f_*c_T}{\epsilon^T(N_{G/Y})}$$

{\bf Proof of the Corollary:} (See also Lemma 2.1 in \cite{LLY2}.) The free
part of $c_T$ is a push-forward of classes from the fixed-point loci of
$X$, by the theorem. Since
$F$ is the only such locus to map to $G$, only its contribution survives
under $j^*f_*c_T$ and under $i^*c_T$, so we may as
well assume that $c_T = i_*b_T$ for $b_T\in H^*(F,{\bf Q})[t,t^{-1}]$. But
then:
$$j^*f_*c_T = j^*f_*i_*b_T = j^*j_*(f|_F)_*b_T =
\epsilon^T(N_{G/Y})(f|_F)_*b_T$$ whereas
$$b_T = \frac{i^*i_*b_T}{\epsilon ^T(N_{F/X})} = \frac{i^*c_T}{\epsilon
^T(N_{F/X})}$$
completing the proof. 

\bigskip

Returning to the basic diagram of the previous section, note that
the ``standard'' linearized action of ${\bf C}^*$ on
${\bf P}^1$:
$$\mu \cdot (a,b) \mapsto (a,\mu b)$$
induces actions on $\overline N_{0,0}({\bf P}^n,d)$ and on
${\bf P}^n_d$ and that
$u:\overline N_{0,0}({\bf P}^n,d) \rightarrow {\bf P}^n_d$
is ${\bf C}^*$-equivariant.
The vector bundles ${\cal O}_{{\bf P}^n_d}(1)$ and $\pi_*e^*E$ are
equipped with natural linearizations, so we may define equivariant chern
classes:
$$h^T := c_1^T({\cal O}_{{\bf P}^n_d}(1)) \ \mbox{and}\ 
[S]_d^T := c_r^T(\pi_*e^*E)$$
on the spaces ${\bf P}^n_d$ and $\overline N_{0,0}({\bf P}^n,d)$
respectively. Then:

\bigskip

\nt {\bf Proposition 3.1:}

\medskip

(1a) $\overline M_{0,1}({\bf P}^n,d) \subset
\overline N_{0,0}({\bf P}^n,d)$ is a component of the fixed-point locus.

\medskip

(1b) $[S]_d$ and $e^*(h)$ extend to the equivariant 
classes
$[S]_d^T$ and $u^*(h^T)$.

\medskip

(2) $\overline M_{0,1}({\bf P}^n,d)$ is the only fixed
component to map to ${\bf P}^n$.

\medskip 

(3) The equivariant Euler classes are:
$$\epsilon^T(N_{\overline
M_{0,1}/\overline N_{0,0}}) = t(t - \psi)\
\mbox{and}\ \epsilon^T(N_{{\bf P}^n/ {\bf P}^n_d}) = \prod_{k=1}^d(h +
kt)^{n+1}$$

{\bf Proof:} (1a,b) and (2) are immediate.
(For a careful treatment, see \cite{BDPP}.) As for
(3), there are inclusions:
$$\overline M_{0,1}({\bf P}^n,d) \subset \overline M_{0,1}({\bf
P}^n,d) \times {\bf P}^1 \subset \overline N_{0,0}({\bf P}^n,d)$$
The (equivariant!) first chern class of the normal bundle to the 
first inclusion is clearly $t$, and the second is $c_1(T_{{\bf P}^1})
- \psi$, which restricts to $t - \psi$ on the fixed-point locus $\overline
M_{0,1}({\bf P}^n,d)$. This gives one Euler class computation. The second
Euler class computation follows from the Euler sequences for the tangent
bundles to ${\bf P}^n$ and ${\bf P}^n_d$ (see \S 7 for a generalization).

\medskip

The correspondence of residues immediately(!) therefore gives us:
$$e_*\left(
\frac{[S]_d}
{t(t - \psi)}\right) = \frac {j^*u_*[S]_d^T}
{\prod_{k=1}^d(h + kt)^{n+1}}$$

Notice that
we get {\bf every} one-point Gromov-Witten invariant
associated to $S$ in this way by expanding $t(t-\psi) = t^{-2} +
t^{-3}\psi + t^{-4}\psi^2 + ...$ so that:
$$\int_{\overline M_{0,1}({\bf P}^n,d)} 
\psi^a \cup e^*(h^b)\cup [S]_d = \ \mbox{coeff of
$t^{-2-a}$ in} \int_{{\bf P}^n} \frac {h^b \cup
j^*{u}_*[S]_d^T} {\prod_{k=1}^d(h + kt)^{n+1}}$$

Thus our challenge is the computation of:
$$j^*{u}_*[S]_d^T \in H^*({\bf P}^n,{\bf Q}) \otimes {\bf
Q}[t]$$

In one case, this is easy. Since $j^*u_*1 = 1$, we have: 

\medskip

\nt {\bf One-Point Invariants of Projective Space:}
$$e_*\left(\frac 1{t(t-\psi)}\right) = 
\frac 1{\prod_{k=1}^d(h + kt)^{n+1}}$$

\medskip

\nt {\bf 4. The Fano Cases:} It is difficult to see how to compute
$j^*u_*[S]_d^T$ when $E$ is a non-split bundle. In this section we will
develop techniques for making the computation when $E$ is a direct sum of
line bundles, and use these techniques to prove the Fano cases of the
mirror conjecture.

\medskip

Let $D \subset \overline N_{0,0}({\bf P}^n,d)$ be the exceptional
divisor for $u: \overline N_{0,0}({\bf P}^n,d)
\rightarrow {\bf P}^n_d$. This is a divisor with normal crossings which
we will eventually describe in detail, but first observe that
there is a birational map
$$\overline N_{0,1}({\bf P}^n,d) \rightarrow {\bf P}^1 \times 
\overline N_{0,0}({\bf P}^n,d)$$ 
and that the exceptional divisor for this map lies over $D$.
Then:

\medskip

\nt {\bf Proposition 4.1:} There is an
equivariant map:
$$\Phi: \pi_*e^*({\cal O}_{{\bf P}^n}(l)) \rightarrow \
\mbox{Sym}^{dl}W^* \otimes u^*{\cal O}_{{\bf P}^n_d}(l)$$
of vector bundles on $\overline N_{0,0}({\bf P}^n,d)$ which is an
isomorphism when restricted to the complement of the boundary divisor
$D$.

\medskip

{\bf Proof:} For $W \cong {\bf C}^2$ and $V \cong {\bf C}^{n+1}$, the
map $v$ below:
$$\overline N_{0,1}({\bf P}^n,d) \stackrel v\rightarrow {\bf P}(W)\times
{\bf P}(\mbox{Sym}^dW^* \otimes V) - -> {\bf P}^n$$
resolves a rational map which is linear in the second factor and has
degree $d$ in the first. The composition is the evaluation map $e$. 

\medskip

\nt {\bf Corollary 4.2:} Let $f:D \hookrightarrow \overline N_{0,0}({\bf
P}^n,d)$ be the inclusion. Then:
$${u}_*[S]_d^T = \prod_{i=1}^m\prod_{k=0}^{dl_i}(l_ih^T + kt) + 
{u}_*f_*c$$ 
for some equivariant class $c$ supported on $D$.

\medskip

{\bf Proof:} This follows from the projection formula, the computation
of the top equivariant chern class of $\mbox{Sym}^{dl}W^*\otimes {\cal
O}_{{\bf P}^n_d}(l)$ and an application of the proposition to the sum of 
vector bundles $\pi_*e^*{\cal O}_{{\bf P}^n}(l_i)$.

\bigskip

\nt {\bf The Fano of Index 2 or More Case:} If 
$l_1 + ... + l_m < n$ then:
$$e_*\left( \frac{[S]_d} {t(t - \psi)}\right) = \frac{j^*u_*[S]_d^T}
{\prod_{k=1}^d(h + kt)^{n+1}} =
\frac{\prod_{i=1}^m\prod_{k=0}^{dl_i} (l_ih +
 kt)} {\prod_{k=1}^d(h + kt)^{n+1}} =: \phi_d$$
i.e. the difference supported on $D$ contributes nothing.

\bigskip

To prove this, we need to understand $D$ better, as
well as the behavior of the  map $\Phi$ when restricted to $D$. The 
self-similar properties of both will lead to this formula and all 
the others. 
Our description will roughly follow Fulton-Pandharipande \cite{FP}:

\medskip

\nt {\bf Components of the boundary divisor:} There
are natural maps:
$$f_i: \widetilde D_i := \overline N_{0,1}({\bf P}^n,d-i) \times_{{\bf
P}^n} 
\overline M_{0,1}({\bf P}^n,i) \rightarrow \overline N_{0,0}({\bf
P}^n,d)$$ 
that desingularize the $d$ components of $D$. The map
$u$ pulls back as follows: 

\medskip

$$\begin{array}{ccc} \overline N_{0,0}({\bf P}^n,d) & 
\stackrel {u}\rightarrow & 
{\bf P}^n_d \\ \\ f_i\uparrow \ \ && \ \ j_i\uparrow \ \ \ \ \ \ \\ \\
\overline N_{0,1}({\bf P}^n,d-i) \times_{{\bf P}^n} 
\overline M_{0,1}({\bf P}^n,i) & \stackrel{\pi_i}\rightarrow & 
{\bf P}^1 \times {\bf P}^n_{d-i} \\ \end{array}$$

Here $\pi_i$ is the 
projection $\rho_i: \widetilde D_i \rightarrow
\overline N_{0,1}({\bf P}^n,d-i)$  followed by the map 
$v$ defined in Proposition 4.1, 
and $j_i$ is the ``multi-linear'' map:
$${\bf P}(W) \times {\bf P}(\mbox{Sym}^{d-i}W^*\otimes V)
\rightarrow {\bf P}(\mbox{Sym}^dW^*\otimes V)$$
(we identify ${\bf P}(W)$ with ${\bf P}(W^*)$
via the canonical $W \cong W^*\otimes \wedge^2W$).

\medskip

\nt {\bf A stratification of the graph space:} We partially stratify
the boundary of
$\overline N_{0,0}({\bf P}^n,d)$ in terms of ``comb'' types, i.e.
sequences $\mu$ of the form:
$$0 \le d_1 < d_2 < ... < d_r < d_{r+1} = d$$

We define $\Delta_i := d_{i+1} - d_{i}$ and
$$\widetilde D_\mu := \overline N_{0,r}({\bf P}^n,d_1)
\times_{({\bf P}^n)^{r}}
\prod_{i=1}^{r} \overline M_{0,1}({\bf P}^n,\Delta_i)$$ 
and note that the natural finite map
$f_\mu:\widetilde D_\mu \rightarrow D_\mu \subset 
\overline N_{0,0}({\bf P}^n,d)$
is a ``desingularization'' of its image $D_\mu$. The general point of
$D_\mu$ is a stable map of degree $d_1$ on the parametrized
component and 
$\Delta_i$ on $r$ other components (the ``teeth'' of 
the comb) each of which meets the parametrized component. Different comb
types with the same $r$ and $d_1$ and the same sets of
$\Delta_i$'s will share the same image, and the map:
$$\coprod_{\mu} \widetilde D_\mu \rightarrow D_\mu$$
from the union over all comb types with image $D_\mu$
has degree $r!$, the order of the permutation group
of the set of $\Delta_i$'s.
The map $u$ pulls back to $\widetilde D_\mu$ according to:
$$\begin{array}{ccc} \overline N_{0,0}({\bf P}^n,d) & 
\stackrel {u}\rightarrow & 
{\bf P}^n_d \\ \\ f_{\mu} \uparrow \ \ && \ \ j_\mu\uparrow 
\ \ \ \ \ \ \\ \\
\widetilde D_\mu &
\stackrel{\pi_\mu}\rightarrow &  ({\bf P}^1)^{r} \times {\bf
P}^n_{d_1} \\
\end{array}$$

Here $\rho_\mu:\widetilde D_\mu \rightarrow \overline
N_{0,r}({\bf P}^n,d_1)$ and $\pi_\mu$ and
$j_\mu$ are defined as before, with the obvious generalized map
$v:\overline N_{0,r}({\bf P}^n,d)
\rightarrow ({\bf P}^1)^{r} \times {\bf P}^n_{d}$. 

\medskip

The following two lemmas are the heart of the Fano proof.

\medskip

\nt {\bf Lemma 4.3:} The fixed-point locus $j:{\bf P}^n \hookrightarrow
{\bf P}^n_d$ factors uniquely through:
$$j:{\bf P}^n \stackrel {j'_\mu}\rightarrow ({\bf
P}^1)^{r} \times {\bf P}^n_{d_1} 
\stackrel{j_\mu}\rightarrow {\bf P}^n_d$$
and if $b$ is an equivariant cohomology class on $({\bf
P}^1)^{r} \times {\bf P}^n_{d_1}$, then:
$$j^*(j_\mu)_*b = (j'_\mu)^*b  \cup
\frac{\prod_{k=1}^d(h + kt)^{n+1}} {t^{r}\prod_{k=
1}^{d_1}
(h + kt)^{n+1}}$$

{\bf Proof:} $j'_\mu$ maps
$(0,...,0) \times {\bf P}^n \hookrightarrow ({\bf P}^1)^{r} 
\times {\bf P}^n
\stackrel {(1,j)}\longrightarrow ({\bf P}^1)^{r} \times {\bf
P}^n_{d_1}$

\medskip

If $j_\mu$ were an embedding, then the second part of the lemma
would be an immediate consequence of the excess intersection formula.
Since it isn't, we will instead use the localization theorem.
The components of the fixed-point locus of the ${\bf C}^*$ action on
${\bf P}^n_d$ are all copies  of ${\bf P}^n$, embedded via Segre as:
$$(e0 + (d-e)\infty) \times {\bf P}^n \subset {\bf P}^d \times {\bf P}^n
\subset  {\bf P}^n_d = {\bf P}(\mbox{Sym}^dW^*\otimes V)$$
and similarly, the components of the fixed-point locus in
$({\bf P}^1)^{r} \times {\bf P}^n_{d_1}$ are all copies
of 
${\bf P}^n$, embedded via:
$$s \times (e0 + (d_1 - e)\infty) \times {\bf P}^n 
\subset ({\bf P}^1)^{r} \times {\bf P}^n_{d_1}$$
where $s\in \{0,\infty\}^r$. Of the
fixed-point loci in this latter space, there is thus only one that 
maps to the image of $j$, namely the image of $j_\mu'$.  By the
localization theorem, an equivariant class on $({\bf P}^1)^{r}
\times {\bf P}^n_{d_1}$ is a sum of its contributions from
fixed-point loci. The contribution  from the image of $j_\mu'$ is
easily computed to be
${j'_\mu}_*b'$ where:
$$b' := \frac{j_\mu '^*
b}{t^{r} 
\prod_{k = 1}^{d_1}(h + kt)^{n+1}}$$

But this is the unique component to map to the image of $j$, hence:
$$j^*{j_\mu}_*b = j^*j_*b' = b' \cup \prod_{k=1}^d(h + kt)^{n+1}$$

\medskip

\nt {\bf Lemma 4.4:} The equivariant virtual class
$[S]_d^T$ decomposes as:
$$[S]_d^T = \sum_\mu \frac 1{r!}{f_\mu}_*(c_\mu \cup
\pi_\mu^*\prod_{i=1}^m\prod_{k=0}^{d_1l_i} (l_ih^T + kt))$$ 
where $h^T$ is pulled back from the projection to ${\bf P}^n_{d_1}$ and
the $c_\mu$ are equivariant cohomology classes which will 
be explicitly described in the proof. 

\medskip

{\bf Proof:} We'll prove this first when $S$ is a hypersurface of degree
$l$. 

\medskip

MacPherson's graph construction \cite{Fu,Mac} gives a decomposition of the
difference 
$[S]^T_d - u^*\prod_{k=0}^{ld} (lh^T + kt)$ as follows.

\medskip

Let $E_d = \pi_*e^*{\cal O}(l)$, $F_d = \mbox{Sym}^{dl}W^*\otimes
u^*{\cal O}_{{\bf P}^n_d}(l)$, and let $G(dl+1,E_d\oplus F_d)$ be the
Grassmann bundle. Consider the locus of scaled graphs:
$$\overline N_{0,0}({\bf P}^n,d) \times {\bf A}^1
\subset G(dl+1,E_d \oplus F_d);\ \ (x,\lambda) \mapsto \Gamma(\lambda
\Phi_x)$$
where $\Gamma(\lambda\Phi_x) = \{(e,\lambda\Phi_x(e)) | e\in {E_d}(x)\}$,
and let $V \subset G(dl+1,E_d\oplus F_d) \times {\bf P}^1$ be the 
closure of this locus. Let $V_\infty \subset G(dl+1,E_d\oplus F_d)$ be the
fiber of $V$ over $\{\infty\} \in {\bf P}^1$. One of the components
of $V_\infty$ is a (reduced) copy of $\overline N_{0,0}({\bf P}^n,d)$
itself, embedded in the Grassmann bundle via the fibers of $F_d$, and all
other components map to proper subvarieties of 
$\overline N_{0,0}({\bf P}^n,d)$. If $V_Z$ is such a component, surjecting
onto $Z \subset \overline N_{0,0}({\bf P}^n,d)$ via the projection
map $\pi_{V_Z}:V_Z \rightarrow \overline N_{0,0}({\bf P}^n,d)$, 
let $n_{V_Z}$ be its multiplicity in $V_\infty$. Then:
$$[S]^T_d - u^*\prod_{k=0}^{ld} (lh^T + kt) = \sum_{V_Z} n_{V_Z}
{\pi_{V_Z}}_*c_{dl+1}(T_{V_Z})$$
where $T_{V_Z}$ is the pull-back of the tautological sub-bundle on 
$G(dl+1,E_d\oplus F_d)$  (see \cite{Fu}, Example 18.1.6).

\medskip

Now notice that $\Phi$ is the push-forward of the inclusion:
$$e^*{\cal O}_{{\bf P}^n}(l) \cong L(-\sum_{i=1}^d ilC_i) \hookrightarrow
L\cong v^*{\cal O}_{{\bf P}^1\times {\bf P}^n_d}(dl,l)$$ 
where $C_i \cong \overline N_{0,1}({\bf P}^n,d-i)\times_{{\bf P}^n}
\overline M_{0,2}({\bf P}^n,i)$
is the unparametrized component of the universal curve over $D_i \subset
D$. 

\medskip

Since the divisors $C_i$ are exactly the exceptional divisors
for the map $\overline N_{0,1}({\bf P}^n,d) \rightarrow {\bf P}^1\times
\overline N_{0,0}({\bf P}^n,d)$ it is clear that the sheaf inclusion holds
with some negative coefficients of the $C_i$. The computation of the
coefficients follows from the observation that the pull-back
$f^*{\cal O}_{{\bf P}^n}(l)$ under a stable map $f:C \rightarrow
{\bf P}^n$ corresponding to a general point of $D_i$ has degree $il$ on the
unparametrized component of $C$.

\medskip

Alternatively, as in the pointwise description of $u$, if we are
given a stable map $f:C \rightarrow {\bf P}^n$ such that the
parametrized component $C_0$ has degree $d_1$ and such that $r$ curves
``bubble'' off this component with degrees $\Delta_1,...,\Delta_r$ at
nodes
$q_1,...,q_r\in C_0$, then $\Phi$ is the following linear map at the point
of $\overline N_{0,0}({\bf P}^n,d)$ corresponding to $f$:
$$H^0(C,f^*{\cal O}(l)) \rightarrow
H^0(C_0,f^*{\cal O}(l)|_{C_0}) \cong H^0({\bf P}^1,{\cal
O}_{{\bf P}^1}(d_1l)) \subset H^0({\bf P}^1,{\cal O}_{{\bf P}^1}(dl))$$ 
where the last inclusion comes from multiplying by the $\Delta_il$th powers
of the equations of the linear forms dual to the $q_i$. It follows that
$\Phi$ enjoys the following pleasant properties:

\medskip

(i) $\Phi$ successively drops rank in codimension one, always along
the self-intersection strata of the boundary divisor $D$. Moreover, if such
a stratum is not a $D_\mu$ (I call these other strata ``hairy combs''),
then it can be embedded in a comb stratum such that
the generic corank of $\Phi$ is the same along the two strata.

\medskip

(ii) $\Phi$ drops rank ``transversally'' along comb-type strata.
That is, if $\Phi$ has generic coranks $m$ and $n$ along two 
comb-type strata which intersect transversally along a third comb-type 
stratum, then $\Phi$ has generic corank $m+n$ along the intersection. 

\medskip

Precisely, $\Phi$ factors as follows when it is pulled back
to $\widetilde D_\mu$:
$$f_\mu^*\Phi :f_\mu^*E_d \rightarrow \rho_\mu^*E_{d_1} 
\stackrel{\rho_\mu^*\Phi}\longrightarrow
\rho_\mu^* F_{d_1} \subset f_\mu^*F_d$$ 
where $E_{d_1},F_{d_1}$ and $\Phi$ on $\overline N_{0,r}({\bf
P}^n,d_1)$ are pulled back from $\overline N_{0,0}({\bf P}^n,d_1)$,
and the map $\rho_\mu^*F_{d_1} \rightarrow f_\mu^*F_d$ is pulled
back via $\pi_\mu^*$ from the projection to 
$({\bf P}^1)^{r} \times {\bf P}^n_{d_1}$  of the sheaf inclusion given by
multiplication by the
$\Delta_il$ powers of the equations of the appropriate diagonals on 
${\bf P}^1 \times ({\bf P}^1)^{r} \times {\bf P}^n_{d_1}$:
$${\cal O}(d_1l,0,...,0,l) \rightarrow
{\cal O}(dl,\Delta_1l,...,\Delta_rl,l)$$

\medskip

For boundary divisors, the kernel of $f_i^*\Phi$ is
pulled back from the kernel of the evaluation map
of bundles on $\overline M_{0,1}({\bf P}^n,i)$:
$$0 \rightarrow E^1_i \rightarrow \pi_*e^*{\cal O}_{{\bf P}^n}(l) 
\rightarrow e^*{\cal O}_{{\bf P}^n}(l) \rightarrow 0$$ 
via the projection $\tau_{i}:\widetilde D_i
\rightarrow
\overline M_{0,1}({\bf P}^n,i)$. If $Z \subset \overline M_{0,1}({\bf
P}^n,i)$ is the image of the section $\sigma_1$, then $E^1_i$ is the first 
in a filtration of sub-bundles:
$$0 = E^{il+1}_i \subset E^{il}_i \subset ... \subset 
E^{k}_i = \pi_*\left(e^*{\cal O}_{{\bf P}^n}(l) \otimes {\cal
O}(-kZ)\right) \subset ... \subset E^1_i \subset E_i$$
which filter ker($f_i^*\Phi$) according to higher-order behavior of $\Phi$
along $D_i$:
$$0 = \tau_i^*E^{il+1}_i \subset \cdots \subset \tau_i^*E^k_i
= f_i^*(\mbox{ker}(\Phi|_{kD_i})) \subset
\cdots \subset \tau_i^*E^1_i \subset f_i^*E_d$$

There is a filtration of $f_i^*F_d$ given by the spans of
the images of $\Phi|_{kD_i}$:
$$\rho_{i}^*F_{d-i} = \pi_i^*F^1_{d-i} \subset \cdots \subset 
\pi_i^*F_{d-i}^{k}
\subset \cdots \subset f_i^*F_d$$
where  $F_{d-i}^k$ is the push-forward of ${\cal O}_{{\bf P}^1\times 
{\bf P}^1 \times {\bf P}^n_{d-i}}((d-i)l + k-1,k-1,l)$.

\medskip

From these descriptions of 
$\Phi$ along boundary divisors (and induction) we can conclude
that infinitessimal versions of the pleasant properties
(i) and (ii) also hold for $\Phi$. Property (i) tells us immediately that
any $Z \subset 
\overline N_{0,0}({\bf P}^n,d)$ in the image of a component $V_Z$ must
be a comb-type boundary stratum, hence that this gives
a decomposition of the form:
$$[S]^T_d = \sum_{D_\mu} \gamma_\mu$$
summed over the images (including $\overline N_{0,0}({\bf P}^n,d)$) of 
the comb-types.
But thanks to the filtrations of $f_i^*E_d$ and $f_i^*F_d$, we
can be much more explicit. Namely, there are $il$ components $V^k_{D_i}$
over each boundary divisor $D_i$, each of which is a birational image of
a ${\bf P}^1$-bundle $V^k_{\widetilde D_i}$ on $\widetilde D_i$. By a local
coordinate computation, 
$V^k_{D_i}$  has multiplicity $k$ in $V_\infty$ and tautological bundle
fitting into an extension of the form:
$$0 \rightarrow \pi_i^*F_{d-i}^k \oplus \tau_i^*E_i^{k+1} \rightarrow
T_{V^k_{\widetilde D_i}}
\rightarrow {\cal O}(-1) \rightarrow 0$$
hence the codimension one contributions $\gamma_i$ have the desired form:
$$\gamma_i = {f_i}_*\left(c_i \cup c_{(d-i)l+1}^T(\pi_i^*F^1_{d-i})\right)
=  {f_i}_*\left(c_i \cup \pi_i^*\prod_{k = 0}^{(d-i)l}(lh^T + kt)\right)$$
where $c_i = \sum_{k=1}^{il}(-k)\tau_i^*c_{il - k}^T(E^{k+1}_i)
\cup \pi_i^*c_{k-1}^T(F_{d-i}^k/F_{d-i}^1)$.

\medskip

The second property of $\Phi$ allows us to construct the
$\gamma_\mu$ inductively. Namely, if we fix a $V^k_{D_i}$
(or rather the corresponding ${\bf P}^1$-bundle over $\widetilde D_i$) and 
apply the MacPherson construction to the map:
$$f_i^*E_d/\tau_i^*E^1_i \cong \rho_i^*E_{d-i}
\stackrel{\rho_i^*\Phi}\rightarrow 
\rho_i^*F_{d-i} = \pi_i^*F_{d-i}^1\left( \subset T_{V^k_{\widetilde
D_i}}\right)$$    then the components of this $V_\infty$ over  
comb type strata of $\overline N_{0,1}({\bf P}^n,d-i)$ map birationally to
the components of the original $V_\infty$ over comb type strata of
$\overline N_{0,0}({\bf P}^n,d)$. It follows by induction that
for each $\mu$ (which includes an ordering of the $\Delta_i$), we will
obtain towers of ${\bf P}^1$-bundles 
mapping birationally to the components $V_{D_\mu}$ over $D_\mu$,
and we obtain the Lemma with the following explicit formula for the
$c_\mu$:
$$c_\mu = \prod_{i=1}^{r} \sum_{k_i = 1}^{\Delta_il}
(-k_i)\tau_{\Delta_i}^*c_{\Delta_il - k_i}^T(E^{k_i + 1}_{\Delta_i})
\cup \pi_\mu^*c_{k_i - 1}^T(F^{k_i}_{d_i}/F^1_{d_i})$$ 
where 
$F^1_{d_i} \subset F^{k_i}_{d_i}$ (we are abusing notation slightly in
the interest of clarity) is pushed forward from the inclusion
of sheaves on ${\bf P}^1\times ({\bf P}^1)^r\times {\bf P}^n_{d_1}$:
$${\cal O}(d_il,\Delta_1l,...,\Delta_{i-1}l,0,...,0,l)
\hookrightarrow {\cal O}(d_il+k_i-1,\Delta_1l,...,\Delta_{i-1}l,
k_i-1,0,...,0,l)$$
and $\tau_{\Delta_i}$ is induced from the projection (so again we are
abusing notation)
$\overline N_{0,1}({\bf P}^n,d_i)
\times_{{\bf P}^n}\overline M_{0,1}({\bf P}^n,\Delta_i) \rightarrow
\overline M_{0,1}({\bf P}^n,\Delta_i)$.

\medskip

In the complete intersection case, the components
$V^k_{\widetilde D_i}$ are only generically ${\bf P}^1$-bundles, 
and the explicit form of the $c_i$ is thus more difficult to 
compute, though the existence of the decomposition of the Lemma
and the inductive nature of the $c_\mu$ terms follow from the same 
argument. Alternatively, one can use the product:
$$[S]^T_d = \prod_{i=1}^m [S_{l_i}]^T_d$$
and the decompositions of each of the hypersurface virtual
classes $[S_{l_i}]^T_d$ along with excess intersection theory to obtain the
desired decomposition of
$[S]^T_d$. There is one subtlety, in that contributions from hairy
combs need to be pushed forward to their
underlying (bald) comb strata.

\nt {\bf Proof of the Fano of Index
$\ge 2$ Case:} Using the decomposition of $[S]_d^T$ from Lemma 4.4, we
have:
$$u_*[S]_d^T = \sum_\mu \frac 1{r!}u_*{f_\mu}_*(c_\mu \cup
\pi_\mu^*b_\mu) = \sum_\mu
\frac 1{r!}{j_\mu}_*\left({\pi_\mu}_*c_\mu \cup
b_\mu\right)$$
where $b_\mu = \prod_{i=1}^m\prod_{k=0}^{d_1l_i}(l_ih^T +
kt)$.
On the other hand, the
codimension of the image of
$j_\mu:({\bf P}^1)^r \times {\bf P}^n_{d_1} \rightarrow {\bf P}^n_d$ is
$$(n+1)(d - d_1) - r \ge n(d - d_1)$$
because of the obvious inequality
$r \le d - d_1$. Thus the codimension of 
${j_\mu}_*b_\mu$ is 
already at least:
$$(\sum_{i=1}^m l_i)d_1 + n(d - d_1) >
\mbox{codim}(u_*[S]_d^T) =  (\sum_{i=1}^m l_i)d$$
since $(\sum_{i=1}^m l_i) < n$ by assumption. This implies 
that ${\pi_\mu}_*c_\mu = 0$ on all strata of positive codimension, 
hence
$u_*[S]_d^T = \prod_{i=1}^m\prod_{k=0}^{dl_i}(l_ih^T + kt)$. Together
with the argument from \S 3, this proves the desired result:
$$e_*\left(\frac{[S]_d}{t(t-\psi)}\right) = \phi_d$$

\nt {\bf The Fano of Index One Case:} Suppose $l_1 + ... + l_m 
= n$. Then:
$$e_*\left(\frac{[S]_d}{t(t-\psi)}\right) = 
\sum_{r=0}^d \frac{(-\prod l_i!)^{r}\phi_{d-r}}{r!t^r}$$

{\bf Proof:} By Lemmas 4.3 and 4.4, the contribution of the 
stratum indexed by $\mu$ to the {\bf ratio}:
$$\frac{j^*u_*[S]_d^T}{\prod_{k=1}^d(h + kt)^{n+1}}$$
is given by: 
$$\frac{j^*{j_\mu}_*({\pi_\mu}_*c_\mu \cup
b_\mu)}{r!\prod_{k=1}^d(h + kt)^{n+1}}
= \frac{j_\mu'^*{\pi_\mu}_*c_\mu \cup
j_\mu'^* b_\mu}{r!t^{r}\prod_{k=1}^{d_1}(h + kt)^{n+1}} 
= \frac{j_\mu'^*{\pi_\mu}_*c_\mu \cup \phi_{d_1}}
{r!t^{r}}$$  
where $b_\mu$ is defined as in the higher index case.

\medskip

When 
$(\sum l_i)d_1 + (n+1)(d - d_1) - r >
(\sum l_i)d$, then as we've already seen, ${\pi_\mu}_*c_\mu = 0$ for
dimension reasons. Under the Fano index one assumption, this inequality
holds unless $r = d - d_1$, in which case equality holds. This
occurs exactly once for  each $r$ from $0$ to $d$, namely when $\mu$
is the comb type:
$$0 \le d_1 < d_1+1 < ... < d_1 + r = d$$

For these comb strata,
${\pi_\mu}_*c_\mu$ is a codimension zero class, so that 
the formula here amounts to the identity:
${\pi_\mu}_*c_\mu =
(-\prod_{i=1}^m l_i!)^{r}\cdot 1.$ 

\medskip

In the hypersurface case, we note that the positive chern classes of
the $\pi_1^*(F^k_{d-1}/F^1_{d-1})$ cannot contribute to the push-forward,
hence the only term which does contribute to ${\pi_1}_*c_1$ is
$-\tau_i^*c_{l-1}(E^2_1)$. The
classes $c_1(E^k_1/E^{k+1}_1)$ are easily computed to be $e^*(lh) +
k\psi$. The $e^*(lh)$ do not contribute to the push-forward, and
we are left with the term $-l!\psi^{l-1}$ which pushes forward to
$-l!$. Using the inductive description of $c_\mu$, we readily obtain
the desired computation ${\pi_\mu}_*c_\mu = (-l!)^r$. 

\medskip

In the complete intersection case, from the product of 
decompositions $[S]^T_d = \prod_{i=1}^m[S_{l_i}]^T_d$
and excess intersection, we can conclude
that modulo terms pulled back under $\pi_1^*$, the term $c_1$ agrees with
$f_1^*(D_1)^{m-1} \prod_{i=1}^m
-\tau_1^*c_{l_i-1}(E^2_{l_i,1})$ where $E_{l_i,1}$ is the bundle 
$\pi_*e^*{\cal O}_{{\bf P}^n}(l_i)$ on $\overline M_{0,1}({\bf P}^n,1)$.
For all $i$, we have $f_i^*(D_i) = -\tau_i^*\psi  - \rho_i^*\psi$, the
second term of which does not contribute to the  push-forward, and
treating the rest of the expression as in the previous paragraph, we
obtain 
$${\pi_1}_*c_1
=
{\pi_1}_*\tau_i^*\left((-\psi)^{m-1}\prod_{i=1}^m(-l_i!\psi^{l_i-1})\right)
= -\prod_{i=1}^m l_i!$$ 

As in the hypersurface case, the desired form of ${\pi_\mu}_*c_\mu$
follows from the inductive description of the $c_\mu$.

\bigskip

\nt {\bf 5. The Calabi-Yau Case.} The decomposition of Lemma 4.4
together with the push-pull formula of Lemma 4.3 always yield:

$$e_*\left(\frac{[S]_d}
{t(t-\psi)}\right) = 
\sum_\mu \frac{\phi_{d_1} \cup {j_\mu'}^*{\pi_\mu}_*c_\mu}{r!t^r}$$
so that the challenge is to find a method for computing the 
${j_\mu'}^*{\pi_\mu}_*c_\mu$. This seems to be very difficult in the 
general-type case, but in the
Calabi-Yau case there is a wonderful simplification:

\medskip
 
\nt {\bf Lemma 5.1:} When $S$ is Calabi-Yau, the decomposition of 
Lemma 4.4 has the following additional properties: 

\medskip

(a) The equivariant chern classes
$\lambda_\mu(h,t) := j_\mu'^*{\pi_\mu}_*c_\mu$ 
satisfy
$$\lambda_\mu(h,t) =
\lambda_{d_1,d_2}(h,t)\lambda_{d_2,d_3}(h,t) \cdots
\lambda_{d_r,d}(h,t)$$
where $\lambda_{d-i,d}(h,t) = j_i'^*{\pi_i}_*c_i$.

\medskip

(b) The ``simple'' classes $\lambda_{d-i,d}(h,t)$ are linear and satisfy:
$$\lambda_{d-i,d}(h,t) = \lambda_{0,i}(h+(d-i)t,t)$$

\medskip

We immediately obtain the following:

\medskip

\nt {\bf A Formula for the Calabi-Yau Case:} If $l_1 + ... + l_m = n+1$,
then there are linear forms
$\lambda_i(h,t)$ only depending upon $(l_1,...,l_m)$ such that:
$$e_*\left(\frac{[S]_d}
{t(t-\psi)}\right) = 
\sum_\mu \frac{\phi_{d_1} \cup
\prod_{i=1}^r \lambda_{\Delta_i}(h + d_it,t)}{r!t^r}$$

\medskip

{\bf Proof:} Set
$\lambda_i(h,t) = \lambda_{0,i}(h,t)$ and apply the lemma.

\bigskip

{\bf Proof of Lemma 5.1:} Recall that Lemma 4.4
gave us:
$$c_\mu = \prod_{i=1}^{r} \sum_{k_i = 1}^{\Delta_il}
(-k_i)\tau_{\Delta_i}^*c_{\Delta_il - k_i}(E^{k_i + 1}_{\Delta_i})
\cup \pi_\mu^*c_{k_i - 1}(F^{k_i}_{d_i}/F^1_{d_i})$$
(again we will start with the hypersurface case) and in
particular,
$$c_i = \sum_{k=1}^{il}(-k)\tau_i^*c_{il - k}(E^{k+1}_i)
\cup \pi_i^*c_{k-1}(F_{d-i}^k/F_{d-i}^1)$$

Recall also the construction of the map $\pi_i$:
$$\pi_i:\overline N_{0,1}({\bf P}^n,d-i)
\times_{{\bf P}^n}\overline M_{0,1}({\bf P}^n,i)
\stackrel{\rho_i}\rightarrow 
\overline N_{0,1}({\bf P}^n,d-i) \stackrel v \rightarrow {\bf P}^1\times
{\bf P}^n_{d-i}$$ 
Let $\xi^T$ be the equivariant hyperplane class on ${\bf P}^1$. Then
one computes:
$$c_{k - 1}^T(F^k_{d-i}/F^1_{d-i}) = 
\prod_{j = 1}^{k-1}(l(h^T + (d-i)t) + j\xi^T)$$

On the other hand, by the K\"unneth decomposition of 
in $\Delta \subset {\bf P}^n\times {\bf P}^n$, the class
${\rho_i}_*\tau_i^*c_{il-k}^T(E^{k+1}_i)$ is of the form
$\kappa_a^i e^*(h^a)$ for a rational number $\kappa_a^i$ that is clearly
independent of $d$. Putting these together, we see that if $S$ is
Calabi-Yau, then ${\pi_i}_*c_i$ is a codimension one class, hence $a=1$ or
$a=0$, and:
$${\pi_i}_*c_i = -v_*\kappa_1^ie^*(h) - 2\kappa_0^i(l(h^T + (d-i)t) +
\xi^T)$$

But $e^*(h) = v^*(h^T + (d-i)\xi^T) - \epsilon$ on $\overline
N_{0,1}({\bf P}^n,d-i)$ for a $v$-exceptional divisor 
$\epsilon$ (from the explicit description of $\Phi$
in Lemma 4.4). Hence:
$${j'_i}^*{\pi_i}_*c_i = -\kappa_1^i(h + (d-i)t) - 2\kappa_0^i(l(h +
(d-i)t) + t)$$
proving part (b) of the lemma with $\lambda_{0,i}(h,t) = (-\kappa_1^i -2
l\kappa_0^i)h -2 \kappa_0^it$.

\medskip

Part (a) is proved similarly. It follows from the projection formula that
$${\rho_\mu}_*c_\mu = 
\prod_{i=1}^{r} \sum_{k_i}
(-k_i){\rho_\mu}_*\tau_{\Delta_i}^*c_{\Delta_il - k_i}^T(E^{k_i +
1}_{\Delta_i})
\cup v^*c_{k_i - 1}^T(F^{k_i}_{d_i}/F^1_{d_i})$$
If we let
$\xi_i^T$ be the hyperplane class of the
$i$th copy of 
${\bf P}^1$, then:
$$c_{k_i-1}^T(F^{k_i}_{d_i}/F^1_{d_i}) = 
\prod_{j=1}^{k_i-1}(l(h^T + d_1t + \Delta_1\xi^T_1 + ... +
\Delta_{i-1}\xi_{i-1}^T) + j\xi_i^T)$$
and as in the proof of (b),
${\rho_\mu}_*\tau_{\Delta_i}^*c_{\Delta_il+1}(E^{2}_{\Delta_i}) = 
\kappa_1^i(v^*(h^T + d_it) - \epsilon_i)$ and 
${\rho_\mu}_*\tau_{\Delta_i}^*c_{\Delta_il+2}(E^{3}_{\Delta_i}) 
= \kappa_0^i$. Different $\epsilon_i$ divisors multiply together to 
yield $v$-exceptional classes, so putting all this together, we get:
$${j_\mu'}^*{\pi_\mu}_*c_\mu = \prod_{i=1}^r
(-\kappa_1^i(h + d_it) - 2\kappa_0^i(l(h + d_it) + t))$$
proving the lemma in the hypersurface case.

\medskip

The complete intersection case is proved similarly. The only 
wrinkle in that case is the presence of normal classes to the strata
in the decomposition. But $f_i^*(D_i) = \tau_i^*\psi + \psi_{{\bf P}^1}
+ \epsilon$ where $\epsilon$ is exceptional for the $v$ map, and thus 
(modulo $\epsilon$) does
not depend upon $d$. Once again, the lemma now follows from the linearity
of the ${\rho_i}_*c_i$ and the inductive description of the $c_\mu$.

\medskip

\nt {\bf 6. Using the Formula.} The real beauty of the 
Calabi-Yau formula:
$$e_*\left(\frac{[S]_d}
{t(t-\psi)}\right) = \sum_\mu \frac{\phi_{d_1} \cup
\prod_{i=1}^r \lambda_{\Delta_i}(h + d_it,t)}{r!t^r}$$
is that it recursively computes itself! This is well-known, but I include
the computation here for the reader's enjoyment. 

\medskip

Namely, consider the simple comb type $\{0 \le 0  < d\}$ which
contributes:
$$\frac{\phi_0\cup \lambda_d(h,t)}t = \frac{\prod_{i=1}^ml_ih \cup
\lambda_d(h,t)}{t}$$ to right side of the formula. This term is
irrelevant to the left side of the formula, which only involves the
powers $t^{-2},t^{-3},...$ so the ``error'' coefficients of $t^{-1}$ and
$t^0$ (there are no more when 
$S$ is Calabi-Yau) in the rest of the right side of the
formula determine this term. Since these only depend upon
$\lambda_1,...,\lambda_{d-1}$, we therefore have an inductive 
construction of the $\lambda_d$. 

\medskip
 
We may read off the coefficients of $\lambda_d(h,t) =
\alpha_dh + \beta_dt$ via:
$$\alpha_dt^{-1} = \frac{1}{\prod_{i=1}^ml_i} \int_{{\bf
P}^n}\frac{h^{n-m-1} \cup \phi_0 \cup \lambda_d(h,t)}{t}$$
$$\beta_d = \frac{1}{\prod_{i=1}^ml_i}\int_{{\bf P}^n}\frac{h^{n-m}\cup
\phi_0 \cup \lambda_d(h,t)}{t}$$

In other words, not only does the formula in degree $d$ compute the 
one-point Gromov-Witten invariants, but via the ``error'' coefficients, it
computes the form $\lambda_d$ which is to be used in higher degrees! 

\medskip

Let $S$ be the quintic threefold in ${\bf P}^4$. Then 
$$e_*\left(\frac{[S]_d}{t(t-\psi)}\right) = n_dh^3t^{-2} +
m_dh^4t^{-3}$$ 
and it follows from the projection
formula that the ``actual'' physical number of rational curves of degree
$d$ on $S$ is $n_d/d = -m_d/2$. For the reader who is unused to this and
nervous about the  fact that these are not typically integers, the
Aspinwall-Morrison  formula (see \cite{AM}) translates these numbers into 
the expected numbers of immersed rational curves
of degree $d$. The formula now produces:

\medskip

$$n_1 = 2875,\ \lambda_1 = - (770)h - (120)t$$
$$n_2 = 4876875/4,\  \lambda_2 = - (421375)h - (60000)t$$
$$n_3 = 8564575000/9,\  \lambda_3 = - (436236875)h -
(59937500)t$$
$$n_4 = 15517926796875/16,\ \lambda_4 = - (17351562078125/6)h -
(390555125000)t$$

Via the Aspinwall-Morrison formula, this translates into:
$$\frac {n_d}d = \sum_{e|d}\frac{N_e}{e^3}$$
where $N_e$ is the expected number of immersed curves of degree $e$. This
gives:
$$N_1 = 2875, N_2 = 609250, N_3 =  317206375\ \mbox{and}\ N_4 =
242467530000$$ 
the well-known numbers of rational curves of degree $\le 4$ on the
quintic.

\medskip

\nt {\bf 7. A Relative Version.} Given a projectivized vector bundle:
$$\pi:{\bf P}(V) \rightarrow  X$$
over a projective manifold $X$, there are relative moduli spaces:
$$\pi_M:\overline M_{0,k}({\bf P}(V),d) \rightarrow X
\ \mbox{and}\ \pi_N:\overline N_{0,k}({\bf P}(V),d) \rightarrow X$$
for stable maps to the fibers of $\pi$. The fibers of 
$\pi_M$ and $\pi_N$ are the moduli spaces we studied previously, and 
these moduli spaces are equipped with evaluation maps (to ${\bf P}(V)$)
and forgetful maps with the usual properties. The definition of 
${\bf P}^n_d$ is readily generalized to:
$${\bf P}(V)_d := {\bf P}(\mbox{Sym}^d(W^*) \otimes V)$$
and the basic diagram in the relative setting is:
$$\begin{array}{ccc} \overline N_{0,0}({\bf P}(V),d)  & \stackrel u
\rightarrow & {\bf P}(V)_d\\ \\
i \uparrow \ \ & & j\uparrow \ \ \\ \\  
 \overline M_{0,1}({\bf P}(V),d) &  \stackrel e \rightarrow &
{\bf P}(V) \end{array}$$
All the results of the paper carry through unchanged, with one exception.
Proposition 3.1(3) now needs to take into account the chern 
classes of $V$. 
Namely, if $\alpha_1,...,\alpha_{n+1}$ are the chern roots of 
$\pi^*V$, then:
$$\epsilon^T(N_{{\bf P}(V)/ {\bf P}(V)_d}) =
\prod_{k=1}^d\prod_{j=1}^{n+1}(h + \alpha_j + kt)$$
which follows from the diagram:
$$\begin{array}{ccccccccc}
& & 0 & & 0 \\
& & \downarrow & & \downarrow \\
& & {\cal O}_{{\bf P}(V)} & = & {\cal O}_{{\bf P}(V)} \\
& & \downarrow & & \downarrow \\
0 & \rightarrow & \pi^*V(1) & \rightarrow & \mbox{Sym}^d(W^*) \otimes
\pi^*V(1)|_{{\bf P}(V)} & 
\rightarrow & N_{{\bf P}(V)/ {\bf P}(V)_d} & \rightarrow & 0\\
& & \downarrow & & \downarrow & & \| \\
0 & \rightarrow & T_{{\bf P}(V)/X} & \rightarrow & 
T_{{\bf P}(V)_d/X}|_{{\bf P}(V)} & \rightarrow  & N_{{\bf P}(V)/ {\bf
P}(V)_d} &
\rightarrow & 0\\
& & \downarrow & & \downarrow \\
& & 0 & & 0 
\end{array}$$ 

We get the following interesting formula already in degree one:

\medskip

\nt {\bf Relative Schubert Calculus:} Chern classes on the relative flag
bundle Fl$(1,2,V)$ over $X$ push forward to ${\bf P}(V)$ via:
$$e_* \left(\frac{\sigma(q_1,q_2)}{t(t-\psi)}\right)
 = \frac{\sigma(h,h+t)} {\prod_{j=1}^{n+1}(h+\alpha_j + t)} +
O(t^{-1})$$ where $\sigma(q_1,q_2)$ is any chern class pulled back from
the relative  Grassmann bundle $G(2,V)$ and expressed as a symmetic
polynomial in the chern roots of the (dual of the) universal
sub-bundle $S^*$. 

\medskip

{\bf Proof:} $S^*$ pulls back to $\pi_*e^*{\cal O}_{{\bf P}(V)}(1)$ on 
Fl$(1,2,V) = \overline M_{0,1}({\bf P}(V),1)$. By Proposition 4.1, the 
map
$\Phi:\pi_*e^*{\cal O}_{{\bf P}(V)}(1) \rightarrow W^*\otimes u^*{\cal
O}_{{\bf P}(V)_d}(1)$ of bundles over $\overline N_{0,0}({\bf P}(V),1)$
is  an isomorphism off the unique boundary stratum $D_1$. It follows
from the argument of \S 3 (and Lemma 4.3) that:
$$e_*\left(\frac{\sigma(q_1,q_2)}{t(t-\psi)}\right) = 
\frac{\sigma(h,h+t)} {\prod_{j=1}^{n+1}(h+\alpha_j + t)} + 
\frac{{j_1'}^*{\pi_1}_*c_1}{t}$$
for some equivariant cohomology class $c_1$ supported on $D_1$.

We can express the denominator in terms of Segre classes:
$$\frac 1{\prod_{j=1}^{n+1}(h+\alpha_j + t)} = 
\frac 1{(h+t)^{n+1}}s_{(h+t)^{-1}}(\pi^*V)$$
where $s_{(h+t)^{-1}}(\pi^*V) = 1 + (h+t)^{-1}s_1(\pi^*V) +
(h+t)^{-2}s_2(\pi^*V) + ...$.

\medskip

In the case of a ``linear'' complete intersection $S \subset {\bf P}(V)$
defined by $m$ sections of ${\cal O}_{{\bf P}(V)}(1)$,
then $[S]_1 = (q_1q_2)^m$
and we get a generalized ``Porteous'' formula for lines:
$$e_*\left(\frac{[S]_1}{t(t-\psi)}\right) =
h^m(h+t)^{m-n-1}s_{(h+t)^{-1}}(\pi^*V) + O(t^{-1})$$

The $t^{-2}$ coefficient is 
$h^m(s_{m-n+1}(\pi^*V) - s_{m-n}(\pi^*V)h)$.
When we multiply this by $h$ and push forward to $X$,
we get the formula:
$$s_{m-n+1}^2(V) - s_{m-n}(V)s_{m-n+2}(V)$$
for the push-forward to $X$ of the class $[S]_1$ on the Grassmann
bundle. This agrees with the classical Porteous formula for lines (see
e.g. \cite{HT}). 

\medskip

For general degree and a general complete intersection $S \subset {\bf
P}(V)$, let:
$$\phi_d =  \frac{\prod_{i=1}^m\prod_{k=0}^{dl_i}(l_ih +
kt)}{\prod_{k=1}^d\prod_{j=1}^{n+1}(h + \alpha_j + kt)}$$

Then the exact analogues of the Fano formulas hold with 
this $\phi_d$, and for Calabi-Yau's we have the analogous:

\medskip

\nt {\bf Relative Calabi-Yau Formula:} If $\sum_{i=1}^m l_i \le n+1$,
then there are linear equivariant classes
$\lambda_e(h,t) \in \ \mbox{H}^*({\bf P}(V),{\bf Q})[t]$ such that:
$$e_*\left(\frac{[S]_d}{t(t-\psi)}\right) = 
\sum_\mu \frac{\phi_{d_1} \cup
\prod_{i=1}^r \lambda_{\Delta_i}(h + d_it,t)}{r!t^r}$$

As an application, consider the ``linear relative Calabi-Yau's'', i.e.
the $S \subset {\bf P}(V)$ that are cut out by $n+1$ transverse sections of
${\cal O}_{{\bf P}(V)}(1)$. The following was first proved in \cite{BT}
in the  context of symmetric products of a smooth curve $C$, where 
the $g-1$st symmetric product $C_{g-1}$ is an example of a linear
relative Calabi-Yau over the Jacobian of $C$:

\medskip

\nt {\bf Proposition 7.1:} The three-point Gromov-Witten invariants:
$$\int_{\overline M_{0,3}({\bf P}(V),d)} e_1^*h \cup e_2^*h \cup 
e_3^*c \cup [S]_d$$
of linear relative Calabi-Yau's are independent of $d \ge 1$. 

\medskip

\nt {\bf Remark:} When $X$ admits no rational curves, this says that
the quantum product of $h$ with itself in the quantum cohomology ring of
$S$ is of the form:
$$h*h = h^2 + bq + bq^2 + bq^3 + ...$$

\nt {\bf Proof:} By the projection formula and the relative Schubert
calculus, the proposition is equivalent to
$e_*([S]_d) = 
(1/d^2)e_*([S]_1) = 
(1/d^2)h^{n+1}(s_2 - s_1h)$.

\medskip

In this case, 
$\phi_d = h^{n+1}s_{(h+t)^{-1}}(\pi^*V)s_{(h+2t)^{-1}}(\pi^*V) ... 
s_{(h+ dt)^{-1}}(\pi^*V)$
$$= h^{n+1}\left(1 + \frac{s_1}t\sum_{k=1}^d \frac 1k +
\frac{s_1^2}{t^2} \sum_{1 \le j < k \le d}\frac 1{jk} + \frac{s_2 -
s_1h}{t^2}\sum_{k=1}^d \frac 1{k^2} + ...\right)$$

\medskip

We make the following (only valid for linear Calabi-Yau's):

\medskip

\nt {\bf Assumption:} The $\lambda_e(h,t)$ of the Calabi-Yau formula are 
independent of $h$.

\medskip

We can separate variables in the formula:
$$e_*\left(\frac{[S]_d}{t(t-\psi)}\right) = 
\sum_{d_1} \phi_{d_1} \sum_{0 < e_1 < ... < e_r = d - d_1}
\frac{\prod_{i=1}^r (\lambda_{e_{i} - e_{i-1}}(t)/t)}{r!}$$
and if we express this in terms of generating functions, we get:
$$h^{n+1} + \sum_{d > 0} q^d e_*\left(\frac{[S]_d}{t(t-\psi)}\right) = 
(\sum_{d \ge 0} q^d\phi_d)\mbox{exp}(\sum_{e > 0}q^e\lambda_e(t)/t)$$

Let $\lambda_e(t) = a_et + b_e$ where $a_e \in {\bf Q}$ and
$b_e$ is a cohomology class of degree one coming from $X$. Then equating
coefficients of $t^0$ gives:
$$1 = (1 + q + q^2 + ...)\mbox{exp}(\sum_e a_eq^e)$$
so that $\sum_e a_eq^e = \ \mbox{log}(1-q)$ and hence $a_e = -\frac 1e$.

\medskip

Similarly, the $t^{-1}$ coefficients give
$0 = \frac{s_1}t \sum_{k=1}^d \frac 1k + \sum_{k=1}^d \frac{b_k}t$
for each $d$, so that by induction, $b_e = -\frac{s_1}e$, and we get
$\lambda_e(t) = -\frac 1e(t + s_1)$ 
which is indeed independent of $h$.
Plug these $\lambda_e(t)$ in the generating function:
$$h^{n+1} + \sum q^de_*[S]_dt^{-2} + ... =$$ 
$$(\sum_e q^e\phi_e) (1 - q)(1 + \frac{s_1}t\mbox{log}(1-q) +
\frac{s_1^2}{2!t^2}\mbox{log}(1-q)^2 + ...)$$
and the $t^{-2}$ term on the right is 
$(1/d^2)h^{n+1}(s_2 - s_1h)$, as desired. 

\medskip

\nt {\bf 8. The Relationship with the Mirror Conjecture.} 
The astute reader will have noticed that the Calabi-Yau formula
in \S 6 does not quite match with the mirror conjecture for Calabi-Yau
complete intersections from the introduction!
We will establish the latter with the aid of the following:

\medskip

\nt {\bf Proposition 8.1:} Let
$$F(q) =  \sum_{d=0}^\infty \sum_{0 < d_1 < ... < d_{r} = d}
\frac{\prod_{i=1}^r (y_{d_i - d_{i-1}} + x_{d_i - d_{i-1}}d_{i-1})}{r!}
q^d$$
Then $\mbox{log}(F(q))$ is a linear function of $y_1,y_2,...$.

\medskip

{\bf Proof:} (Pavel Etinghof) Let $E$ be the Euler vector field:
$$E = \sum_{k=1}^\infty kx_k\frac{\partial}{\partial x_k} + 
ky_k \frac{\partial}{\partial y_k}$$
Then we may rewrite $F(q)$ as:
$$F(q) = \ \mbox{exp}({\sum_{d=1}^\infty q^dy_d + q^dx_dE}) \cdot 1 =
\sum_{r=0}^\infty
\frac{(\sum q^dy_d + q^dx_dE)^r}{r!} \cdot 1$$ 

Consider the function 
$G(t,q) = \ \mbox{exp}(t{\sum_{d=1}^\infty q^dy_d +
q^dx_dE}) \cdot 1$  
satisfying $G(1,q) = F$, \ $G(0,q) = 1$ and
$q \frac{\partial G}{\partial q} = 
E\cdot G$ (because the Euler vector field is homogeneous of degree zero).
Thus:
$$\frac{\partial G}{\partial t} = \sum_{d=1}^\infty q^dy_dG + 
\sum_{d=1}^\infty q^dx_dE\cdot G$$
hence
$$\frac{\partial \mbox{log}(G)}{\partial t} = 
\frac 1G \frac{\partial G}{\partial t} = 
\sum_{d=1}^\infty q^dy_d + \sum_{d=1}^\infty q^{d+1} x_d \frac{\partial
\mbox{log}(G)}{\partial q}.$$

By the fundamental theorem of calculus(!) this gives:
$$\mbox{log}(G)(t,q) = \sum_{d=1}^\infty tq^dy_d + 
\sum_{d=1}^\infty q^{d+1}x_d 
\int_0^t\frac{\partial \mbox{log}(G(s,q))}{\partial
q}\ ds$$ which proves the desired linearity of log$(F) = \
\mbox{log}(G)(1,q)$ in the $y$ variables by induction on the power of $q$.

\bigskip

\nt {\bf Corollary:} If we let
$F(q) =  \ \mbox{exp}(\sum_{d=1}^\infty y'_dq^d)$ then
$$y'_d = \sum_{0 < d_1 < ... < d_r = d}
\frac{ y_{d_1}\prod_{i=2}^r (x_{d_i - d_{i-1}}d_{i-1})}{r!}$$

{\bf Proof:} Cast out the non-linear terms (in the
$y_k$) from the identity:
$$\sum_{0 < d_1 < ... < d_r = d}
\frac{\prod_{i=1}^r y'_{d_i - d_{i-1}}}{r!}
 = \sum_{0 < d_1 < ... < d_r = d}
\frac{\prod_{i=1}^r (y_{d_i - d_{i-1}} + x_{d_i - d_{i-1}}d_{i-1})}{r!}$$

\bigskip

Finally, suppose $S$ is Calabi-Yau and let $\lambda_d = \alpha_dh +
\beta_dt$. Then our formula may be written as follows:
$$\Sigma(q) = \sum_{d}\phi_d\sum_{0 < e_1 < ... < e_r = e}
\frac{\prod_{i=1}^r (\alpha_{e_i -
e_{i-1}}(d + \frac ht) + \beta_{e_i - e_{i-1}}) + \alpha_{e_i -
e_{i-1}}e_{i-1}}{r!}q^{d+e}$$ 

The proposition applies to give us the new coordinates $y'_e$ which are
linear in $\alpha_e(d+\frac ht) + \beta_e$ (and polynomial in the 
$\alpha_k$'s), hence of the form:
$$y'_e = a_e(d + \frac ht) + b_e$$
where $a_e$ and $b_e$ are independent of $d$. This gives:
$$\Sigma(q) = \sum_{d}\phi_d\sum_{0 < e_1 < ... < e_r = e}
\frac{\prod_{i=1}^r (a_{e_i -
e_{i-1}}(d + \frac ht) + b_{e_i - e_{i-1}})}{r!}q^{d+e}$$ 
which proves the mirror conjecture.

\bigskip

\nt University of Utah, Salt Lake City, UT 84112

\medskip

\nt bertram@math.utah.edu
\end{document}